\documentclass[10pt,leqno]{amsart}
\usepackage{amscd}
\usepackage{amssymb}
\usepackage{amsfonts}
\usepackage{latexsym}
\usepackage{verbatim}
\usepackage{amsthm}
\numberwithin{equation}{section}

\theoremstyle{plain}
\newtheorem{theorem}{Theorem}[section]
\newtheorem{definition}[theorem]{Definition}
\newtheorem{lemma}[theorem]{Lemma}
\newtheorem{prop}[theorem]{Proposition}
\newtheorem{cor}[theorem]{Corollary}
\newtheorem{rem}[theorem]{Remark}

\newcommand{\ov}[1]{\overline{#1}}

\newcommand{\de}[2]{\frac{\partial #1}{\partial #2}}

\newcommand{\Om}{\Omega}
\newcommand{\e}{\epsilon}
\begin{document}
\title{The Ricci tensor of SU(3)-manifolds}
\author{Lucio Bedulli and Luigi Vezzoni}
\date{\today}
\address{Dipartimento di Matematica ''L. Tonelli''\\ Universit\`a di Pisa\\
Largo B. Pontecorvo 5\\
56127 Pisa\\ Italy} \email{vezzoni@mail.dm.unipi.it}
\address{Dipartimento di Matematica - Universit\`a di Bologna\\
Piazza di Porta S. Donato 5\\40126 Bologna\\Italy}
\email{bedulli@math.unifi.it}
\thanks{This work was supported by the Project M.I.U.R. ``Geometric Properties of Real and Complex Manifolds'' and by G.N.S.A.G.A.
of I.N.d.A.M.}
\begin{abstract}
Following the approach of Bryant \cite{B2} we study the intrinsic
torsion of a SU$(3)$-manifold deriving a number of formulae for the
Ricci and the scalar curvature in terms of torsion forms. As a consequence 
we prove that in some special cases the Einstein condition forces the vanishing 
of the intrinsic torsion.
\end{abstract}
\maketitle
\newcommand\C{{\mathbb C}}
\newcommand\R{{\mathbb R}}
\newcommand\Z{{\mathbb Z}}
\newcommand\T{{\mathbb T}}
\newcommand\GL{{\rm GL}}
\newcommand\SL{{\rm SL}}
\newcommand\SO{{\rm SO}}
\newcommand\Sp{{\rm Sp}}
\newcommand\U{{\rm U}}
\newcommand\SU{{\rm SU}}
\newcommand{\Gdue}{{\rm G}_2}
\newcommand\re{\,{\rm Re}\,}
\newcommand\im{\,{\rm Im}\,}
\newcommand\id{\,{\rm id}\,}
\newcommand\tr{\,{\rm tr}\,}
\renewcommand\span{\,{\rm span}\,}
\newcommand\Ann{\,{\rm Ann}\,}
\newcommand\Hol{{\rm Hol}}
\newcommand\Ric{{\rm Ric}}

\newcommand\nc{\widetilde{\nabla}}
\renewcommand\d{{\partial}}
\newcommand\dbar{{\bar{\partial}}}
\newcommand\s{{\sigma}}
\newcommand\sd{{\bigstar_2}}
\newcommand\K{\mathbb{K}}
\renewcommand\P{\mathbb{P}}
\newcommand\D{\mathbb{D}}

\newcommand\f{{\varphi}}
\newcommand\g{{\frak{g}}}
\renewcommand\k{{\kappa}}
\renewcommand\l{{\lambda}}
\newcommand\m{{\mu}}
\renewcommand\O{{\Omega}}
\renewcommand\t{{\theta}}
\newcommand\ebar{{\bar{\varepsilon}}}
\section*{Introduction}
In the last years geometric and physical motivations led
many mathematicians to focus on the geometry of $\mbox{SU}(3)$ and
$\mbox{G}_2$-structures on 6 and 7-dimensional manifolds and on the
interplay between them (see e.g. \cite{AF}, \cite{AFNP}, \cite{AFS}, \cite{AS},
\cite{B2}, \cite{CabSw}, \cite{CS} , \cite{CSw}, \cite{Cley}, \cite{FI} 
and the references therein).
New directions in this field were suggested by
the work of Hitchin \cite{Hi2}.
The present work is inspired by
\cite{B2}, where the author computes the Ricci curvature of a
G$_2$-structure
in terms of the derivatives of the defining 3-form.

In this paper we study the intrinsic torsion of SU$(3)$-manifolds
relating it to the curvature of the induced metric.\\
A SU$(3)$-structure on a 6-dimensional manifold is determined by a pair
$(\kappa,\Omega)$, where $\kappa$ is an almost symplectic structure
and $\Omega$
is a normalized $\kappa$-positive $3$-form (see Section \ref{SU3section} for the definition).
In fact such a pair induces a
natural $\kappa$-calibrated
almost complex structure $J$ on $M$ such that the complex valued
form
$$
\varepsilon=\Omega+i\,J\Omega
$$
is of type (3,0) with respect to $J$.
The intrinsic torsion
of a SU(3)-structure can be described in terms of the derivatives of
the defining forms $(\kappa,\Omega)$ by considering a natural
decomposition of $\Lambda^3M$ and $\Lambda^4M$ in irreducible
SU(3)-submodules. Namely the forms
$d\kappa$, $d\Omega$ and $d^*\Omega$ decompose as
$$
\begin{aligned}
d\kappa&=-\frac{3}{2}\sigma_0\,\Omega+\frac{3}{2}\pi_0\,J\Omega+\nu_1\wedge\kappa+\nu_3\,;\\[3pt]
d\Omega& =\pi_0\,\kappa^2+\pi_1\wedge\Omega-\pi_2\wedge\kappa\,;\\[3pt]
dJ\Omega&=\sigma_0\,\kappa^2+J\pi_1\wedge\Omega-\sigma_2\wedge\kappa\,,
\end{aligned}
$$
where
$\pi_0,\sigma_0,\pi_1,\nu_1,\sigma_2,\nu_3$
lie in different SU(3)-modules.
The forms
$\{\pi_0,\sigma_0,\pi_1,\nu_1,\sigma_2,\nu_3\}$
are called the
\emph{torsion forms} and they vanish
if and only if the SU$(3)$-structure is integrable,
i.e. if and only if the induced metric is Ricci-flat so that
$(M,\kappa,\Omega)$ is a Calabi-Yau threefold.
Moreover special non-integrable SU(3)-structures, e.g.
generalized Calabi-Yau structures\footnote{We
remark that the notion of generalized Calabi-Yau structure we
consider is the one adopted in \cite{dBT} which is different from that one given by
Hitchin in \cite{Hi1}.} and half-flat structures, can be
characterized in terms of torsion forms. In the spirit of \cite{B2}
a principal bundle approach allows us to write down the Ricci tensor
and the scalar curvature of a SU(3)-manifold in terms of torsion
forms. As a direct consequence of these formulae we get that the
scalar curvature of a generalized Calabi-Yau manifold is
non-positive and it vanishes identically if and only if the SU(3)-structure is
integrable. We also prove that the metric of a special generalized Calabi-Yau manifold $M$ is 
Einstein if and only if $M$ is a genuine Calabi-Yau manifold.
\newline
The paper is organized as follows. In section \ref{SUnsection}
general SU$(n)$-structures are introduced. In section
\ref{SU3section}, which is the algebraic core of the paper, we
specialize to the 6-dimensional case studying the algebra underlying
SU$(3)$-structures. In particular we exhibit an explicit expression
for the complex structure induced by $(\kappa,\Omega)$. In this
section we define the torsion forms and characterize various special
SU$(3)$-structures in terms of these forms. The work in section
\ref{Riemanniansection} follows the steps of \cite{B2} where the
formula for the Ricci curvature of a G$_2$-structure is derived. We
exploit the algebraic formulae obtained in section \ref{SU3section} in order to
come to the explicit formula for the Ricci tensor
\eqref{riccitensor}. Here the final computation was carried out with
the aid of {\sc Maple} while a representation-theoretic argument
justifies the final formulae. In section \ref{SGCYsection} we
collect the above mentioned consequences of formula \eqref{riccitensor} in the 
special case of generalized Calabi-Yau manifolds.  
Section \ref{examplesection} is devoted to the
explicit computations performed on a non-integrable special
generalized Calabi-Yau nilmanifold which illustrate the role of the
torsion forms in this case. In the appendix some technical proofs are provided.\\
\newline
{\sc Acknowledgments:} The authors are grateful to Robert Bryant for
supplying them with the computer programs he used to perform the symbolic
computations in the G$_2$-case. 
They are also grateful to Richard Cleyton for suggesting a considerable 
strengthening of a previous version of Corollary \ref{Cleyton}.\\
\newline
{\sc Notations.} Given a manifold $M$, we denote by $\Lambda^rM$ the
space of smooth $r$-forms on $M$ and we set
$\Lambda^{\bullet}M:=\bigoplus_{r=1}^n \Lambda^r M$. When an almost complex structure
$J$ on $M$ is given, $\Lambda_J^{p,q}M$ denotes the space of complex
forms on $M$ of type $(p,q)$ with respect to $J$.\\
The symplectic group, i.e.\ the group of automorphisms of
$\R^{2n}$ preserving the standard symplectic form $\kappa_n=\sum_{i=1}^n
dx_{2i-1}\wedge dx_{2i}$, will be denoted by $\mbox{Sp}(n,\R)$.\\
Furthermore when a coframe $\{\alpha_1,\dots,\alpha_n\}$ is given we
will denote the $r$-form $\alpha_{i_{1}}\wedge\dots\wedge\alpha_{i_{r}}$ by
$\alpha_{i_{1}\dots i_{r}}$.\\
In the indicial expressions the symbol of sum over repeated indices
is omitted.

\section{SU$(n)$-structures}
\label{SUnsection}
\subsection{U$(n)$-structures.}
Let $(M,\kappa)$ be a $2n$-dimensional almost symplectic manifold.
The \emph{symplectic Hodge operator}
$$
\bigstar\colon \Lambda^rM\to\Lambda^{2n-r}M\,,
$$
is defined by means of the relation
$$
\alpha\wedge\bigstar\beta=\kappa(\alpha,\beta)\frac{\kappa^n}{n!}\,,
$$
where $\alpha,\beta\in \Lambda^rM$. It is easy to check that
$\bigstar^2=I$.
An almost complex structure on $M$ is an endomorphism $J$ of $TM$ such that
$J^2=-I$. Note that the endomorphism induced by $J$ on
$\Lambda^pM$ (again denoted by $J$)  satisfies the identity $J^2=(-1)^pI$.
An almost complex structure is said to be $\kappa$-\emph{tamed} if
$$
\kappa_x(v,J_xv)> 0
$$
for every $x\in M$ and non-zero vector $v\in T_xM$.
If further $\kappa$ is preserved by $J$, the almost complex structure is said
to be $\kappa$-\emph{calibrated}. In this case we denote by $g_J$ the
Riemannian metric
\begin{equation}
\label{g}
g_{J}(X,Y):=\kappa(X,JY)\,,
\end{equation}
for every vector field $X,Y$ on $M$.
We immediately get that $J$ is an isometry of $g_J$, i.e. $g_J$ is $J$-Hermitian.
We denote by $\mathcal{C}_{\kappa}(M)$ the space of
$\kappa$-calibrated almost complex structures on $M$.
The elements of $\mathcal{C}_{\kappa}(M)$ can be viewed as smooth global
sections of a fiber bundle whose fibers are isomorphic to the homogeneous space
$$
\mbox{Sp}(n,\R)/\mbox{U}(n)
$$
(see e.g. \cite{AL}).
Since the latter is topologically a $(n+n^2)$-dimensional cell,
given any almost symplectic form $\kappa$, there are always plenty
of $\kappa$-calibrated almost complex structures. Furthermore the
fact that $\mathcal{C}_{\kappa}(M)$ is contractible makes it
possible to define the first Chern class $c_1(M,\kappa)$ of the almost
symplectic manifold $(M,\kappa)$ as $c_1(M,J)$, where
$J \in \mathcal{C}_{\kappa}(M)$.

Given $J\in\mathcal{C}_{\kappa}(M)$ the complexified exterior algebra
$\Lambda^{\bullet} M\otimes\C$ is $\Z^+$-bigraded with respect
to the type as
$$
\Lambda^{\bullet} M\otimes\C=\bigoplus_{r=0}^{2n}\bigoplus_{p+q=r}\Lambda^{p,q}_J M\,.
$$
The metric $g_J$ together with the orientation given by $\kappa$
defines also the classical \emph{Hodge operator},
that in this setting is a $\C$-linear map
$*\colon\Lambda_J^{p,q}M\to\Lambda_J^{n-q,n-p}M$, such that
\[
\alpha\wedge\overline{*\beta}=g_J(\alpha,\overline{\beta})\frac{\kappa^n}{n!}\,,
\]
for all $\alpha,\beta\in \Lambda^{p,q}_JM$. It is well known that $*$
commutes with $J$ and that their composition equals the $\C$-linear
extension of the symplectic Hodge operator:
$$
*J=J*=\bigstar\,.
$$
Since we have
\[
d\colon\Lambda_J^{p,q}M \to \Lambda_J^{p+2,q-1}M \oplus \Lambda_J^{p+1,q}M
\oplus \Lambda_J^{p,q+1}M \oplus \Lambda_J^{p-1,q+2}M\,,
\]
the exterior differential operator accordingly splits as
\[
d=A_J+\partial_J+\bar{\partial}_J+\bar{A}_J.
\]
It is well known that an almost complex structure is integrable if and
only if $\bar{A}_{J}=0$.
 \subsection{SU$(n)$-structures.}
Let $M$ be a $2n$-dimensional manifold and $\mathcal{L}(M)$ be the
GL$(2n,\R)$-principle bundle of linear frames. A SU$(n)$-structure on
$M$  is a SU$(n)$-reduction of $\mathcal{L}(M)$.
Since SU$(n)$ is the group of the unitary transformation of $\C^n$
preserving the standard complex volume form, a SU$(n)$-structure on $M$ is
determined by the choice of the following data:
\begin{itemize}
\item an almost complex structure $J$ on $TM$;
\vspace{0.1 cm}
\item a $J$-Hermitian metric $g$;
\vspace{0.1 cm}
\item a complex $(n,0)$-form $\varepsilon$ of constant
norm $2^\frac{n}{2}$.
\end{itemize}
Alternatively these data can be replaced by
\begin{itemize}
\item an almost symplectic structure $\kappa$;
\vspace{0.1 cm}
\item a $\kappa$-calibrated almost complex structure $J$;
\vspace{0.1 cm}
\item a complex $(n,0)$-form $\varepsilon$, satisfying $\varepsilon \wedge
\overline{\varepsilon} = c_n\,\frac{\kappa^n}{n!}$,
with $c_n=(-1)^\frac{n(n+1)}{2}(2i)^n$;
\end{itemize}
where $\kappa$ and $g$ are relied by  \eqref{g}.
Denote by $\nabla$ the Levi-Civita connection induced by
$g$ on $TM$. We will say that a SU$(n)$-structure is \emph{integrable} if
the restricted holonomy group
Hol$^{0}(TM,\nabla)$ is isomorphic to a subgroup of SU$(n)$.\\
Since the holonomy is determined by the parallel tensors, a
SU$(n)$-structure is integrable if the corresponding triple
$(\kappa,J,\varepsilon)$ satisfies
$$
\nabla\kappa=0\,,\quad\nabla J=0\,,\quad\nabla\varepsilon=0\,.
$$
In this case $(M,\kappa,J,\varepsilon)$ is said to be a \emph{Calabi-Yau
  manifold}.\\
\begin{rem}\label{rem1}
\emph{
Let $(M,\kappa,J,\varepsilon)$ be a SU$(n)$-manifold and assume
\[
d\kappa=0\,,\quad
d\varepsilon=0\,,
\]
then $(M,\kappa,J,\varepsilon)$ is a
Calabi-Yau manifold. In fact if $\alpha\in\Lambda^{1,0}_JM$ we have
$$
0=d(\varepsilon\wedge\alpha)=(-1)^n \varepsilon\wedge d\alpha=(-1)^n
\varepsilon\wedge \overline{A}_J\alpha\,,
$$
hence $\ov{A}_J=0$, which implies that $J$ is integrable. Furthermore, since $\kappa$
is closed, the pair $(\kappa,J)$ defines a K\"ahler structure on $M$;
hence we get
$$
\nabla\kappa=0\,,\quad \nabla J=0\,.
$$
Finally the equation
$\varepsilon\wedge\overline{\varepsilon}=c_n\,\frac{\kappa^n}{n!}$ forces
$\varepsilon$ to be parallel.
}
 \end{rem}
Several non-integrable SU$(n)$-structures are worth to be considered
for both geometrical and physical reasons (the survey article \cite{A} is a good
reference for recent results on non-integrable geometries).\\
A notion of generalized Calabi-Yau manifold has been introduced by
de Bartolomeis and Tomassini; in \cite{dBT} they give the following
definition:
\begin{definition}
\label{GCY}
A \emph{generalized Calabi-Yau (GCY)   structure} on $M$ is a \emph{SU}$(n)$-structure
$(\kappa,J,\varepsilon)$ satisfying the following conditions:
\begin{enumerate}
\item[1.]$d\kappa=0$ $($i.e. $(M,\kappa)$ is a symplectic manifold$)$;
\vspace{0.1 cm}
\item[2.]$\overline{\partial}_J\varepsilon=0$.
\end{enumerate}
\end{definition}
We emphasise again that a different generalization of Calabi-Yau structures
has been considered by Hitchin in a broader context in \cite{Hi1}.
\begin{rem}\emph{
For an almost K\"ahler manifold (i.e. a symplectic manifold endowed
with a calibrated almost complex structure)
it is natural to consider on $TM$  the canonical
Hermitian connection $\widetilde{\nabla}$, whose covariant
derivative is given by
$$
\widetilde{\nabla}_X=\nabla_X-\frac{1}{2}J\nabla_X J\,.
$$
It is characterized by the following properties
$$
\widetilde{\nabla}\kappa=0\,,\quad\widetilde{\nabla} J=0\,,
\quad T^{\widetilde{\nabla}}=\frac{1}{2}\,N_J\,,
$$
where $N_J$ is the Nijenhuis tensor associated to $J$ and
$T^{\widetilde{\nabla}}$ is the torsion of $\widetilde{\nabla}$.
This connection coincides
with $\nabla$ if and only if the pair
 $(\kappa,J)$ is a K\"ahler structure on $M$
(i.e. if and only if $J$ is integrable).}

\emph{If $(M,\kappa,J,\varepsilon)$ is a symplectic
  SU$(3)$-manifold, then the constraint
$\varepsilon \wedge \overline{\varepsilon}=c_n\frac{\kappa^n}{n!}$ implies
$$
\overline{\partial}_J\varepsilon=0\
\iff \widetilde{\nabla}\varepsilon=0\,,
$$
(see \cite{dBT}). Hence GCY manifolds can be defined as ${\rm SU}(n)$-manifolds
with the volume
form $\varepsilon$ satisfying $\widetilde{\nabla}\varepsilon=0$. It
follows that in the GCY case the holonomy group
$\mbox{Hol}^0(TM,\widetilde{\nabla})$ is isomorphic to a subgroup of
SU$(n)$.}
\end{rem}
\section{SU(3)-structures}
\label{SU3section}
In this section we specialize to the case $n=3$ and study the linear
algebra underlying SU$(3)$-structures. Fix a real
6-dimensional symplectic vector space $(V,\kappa)$. Let us denote
by Sp$(V,\kappa)$ the group of automorphisms of the pair
$(V,\kappa)$, i.e. Sp$(V,\kappa)=\{\phi\in \mbox{GL}(V)\,\colon\, \phi^*\kappa=\kappa\}$.
The space of skew-symmetric 3-forms on $V$ splits into the following
two irreducible Sp$(V,\kappa)$-modules
$$
\begin{aligned}
&\Lambda^{3}_0V^*=\{\phi\in\Lambda^3V^*\,|\,\phi\wedge\kappa=0\}\,,\\
&\Lambda^{3}_6V^*=\{\alpha\wedge\kappa\,|\,\alpha\in V^*\}\,.
\end{aligned}
$$
The 3-forms lying in the space  $\Lambda^3_0V^*$ are sometimes
called in the literature \emph{effective} 3-forms (see e.g. \cite{Banos}).
Let us consider the action $\Theta$ of the Lie group
$G=$Sp$(V,\kappa)\times \R^*_+$ on the space $\Lambda^{3}_{0} V^*$
given by
$$
\Theta(\phi,t)\cdot\alpha:=t\,(\phi^{-1})^*\alpha\,,
$$
where $\R^*_+$ denotes the group of positive real numbers.
It is known that this action has an open orbit
$\mathcal{O}$ whose isotropy is locally isomorphic to SU(3) (see
e.g. \cite{Banos} and \cite{Sato}). We will call $\kappa$-\emph{positive} 3-forms
the elements of the orbit $\mathcal{O}$. Since the stabilizer at $\Omega\in\mathcal{O}$ is
locally isomorphic to SU(3), each $\kappa$-positive 3-form singles out a
$\kappa$-calibrated complex structure on $V$ which we are able to
explicitly write down. In fact we have:
\begin{prop}
The endomorphism $P_{\Omega}$ of $V^*$ given by
$$
P_\Omega\colon \alpha\longmapsto
-\frac{1}{2}\bigstar(\Omega\wedge\bigstar(\Omega\wedge \alpha))
$$
has the following properties
\begin{enumerate}
\item[1.] $P_{\Omega}^2$ is a negative multiple of the identity;
\vspace{0.1 cm}
\item[2.] $\kappa(P_{\Omega}\,\alpha,\beta)=-\kappa(\alpha,P_{\Omega}\beta)$, for every $\alpha,\beta\in\Lambda^1V^*$.
\end{enumerate}
\end{prop}
\begin{proof}
1. First we observe that $P_{\Omega}$ is a SU$(3)$-invariant
endomorphism of $V^*$, since it is built using only $\Omega$ and
$\bigstar$. Since SU$(3)$ acts irreducibly on $V^*$, the real
version of Schur's lemma assures that $P_{\Omega}=a\,I+b\,J$, where $J$
is a complex structure on
$V^*$ and $a,b$ are real numbers. \\
Now we claim that $P_{\Omega}^2$ has a negative eigenvalue. From this
claim the conclusion follows. Suppose indeed that there exists $v\neq
0$ such that $P_{\Omega}^2v=\lambda\,v$, with $\lambda<0$. Then
$$
2ab\,Jv=(\lambda^2-a^2+b^2)\,v\,.
$$
If $ab\neq 0$, then $J$ would have a real eigenvalue and this is
impossible. On the other hand if $b=0$ then $P_{\Omega}^2=a^2I$, which
is a contradiction with the claim. Hence $P_{\Omega}=bJ$.
To prove the claim we must use an explicit frame $\{e^1,\dots,e^{6}\}$
of $V^*$ in which $\kappa$ and $\Omega$ takes the standard form
and perform the computation e.g. of $P_{\Omega}^2e^1$.\\

2. We have
\begin{equation*}
\begin{split}
\kappa(P_{\Omega}\alpha,\,\beta)\frac{\kappa^{3}}{6}=&-\kappa(\beta,\,P_{\Omega}\alpha)\frac{\kappa^{3}}{6}
=\frac{1}{2}\beta\wedge\Omega\wedge\bigstar(\Omega\wedge \alpha)=\\
=&-\frac{1}{2}\kappa(\beta\wedge\Omega,\,\alpha\wedge\Omega)\frac{\kappa^{3}}{6}
=-\frac{1}{2}\kappa(\alpha\wedge\Omega,\,\beta\wedge\Omega)\frac{\kappa^{3}}{6}=\\
=&\kappa(P_{\Omega}\beta,\,\alpha)\frac{\kappa^{3}}{6}=-\kappa(\alpha,\,P_{\Omega}\beta)\frac{\kappa^{3}}{6}\,.
\end{split}
\end{equation*}
\end{proof}
It follows:
\begin{cor}
The endomorphism $J_{\Omega}$ $\kappa$-dual to
$($\emph{det}$P_{\Omega})^{-\frac{1}{6}}P_{\Omega}$ is a
$\kappa$-calibrated almost complex structure on $V$.\\
Furthermore the form
$$
\varepsilon=\Omega+iJ_{\Omega}\Omega
$$
is a complex form of type $(3,0)$ with respect to $J_{\Omega}$. If further
$\det(P_\Omega)=1$, then
\begin{equation}
\label{cienne}
\varepsilon\wedge\ov{\varepsilon}=i\frac{4}{3}\,\kappa^{3}\,.
\end{equation}
\end{cor}
We have also this characterization of $\kappa$-positive 3-forms
\begin{lemma}
\label{facile}
These facts are equivalent
\begin{enumerate}
\item[1.] $\Omega$ is a $\kappa$-positive $3$-form;
\vspace{0.1 cm}
\item[2.] the map $F_{\Omega}\colon \Lambda^1V^*\ni\alpha\mapsto\alpha\wedge\Omega$ is
  injective and $\kappa$ is negative definite on the image of $F_{\Omega}$.
\end{enumerate}
\end{lemma}
\begin{rem}
\emph{Note that since $\kappa$ is $J_{\Omega}$-invariant, also
  $J_{\Omega}\Omega$ is effective, i.e. $\kappa\wedge J_{\Omega}\Omega=0$.
}
\end{rem}
\begin{definition}
A $\kappa$-positive 3-form is said to be {\em normalized} if $\det(P_\Omega)=1$.
\end{definition}
\noindent From now on we will drop the subscript $\Omega$ from $J_{\Omega}$
when no confusion arises.\\
\newline
In order to make the exposition more concrete we identify $V$ with
$\R^6$; we denote by $\{e_1,\dots,e_6\}$ the standard basis
and by  $\{e^1,\dots,e^6\}$ the dual one.\\
Fix on $V$ the standard symplectic form
$$
\kappa_0=e^{12}+e^{34}+e^{56}\,
$$
and the standard complex volume form
$$
\varepsilon_0=(e^1+ie^2)\wedge(e^3+ie^4)\wedge(e^5+ie^6)\,.
$$
The real part of $\varepsilon_0$
$$
\Omega_0=e^{135}-e^{146}-e^{245}-e^{236}
$$
is a normalized $\kappa_0$-positive 3-form. The complex structure
associated to $\Omega_0$ is exactly the standard $\kappa_0$-calibrated
complex structure $J_0$ defined by
$$
J_0(e_1)=e_2,\quad J_0(e_3)=e_4,\quad J_0(e_5)=e_6\,.
$$
We will denote by $g_0$ the scalar product associated to
$(\kappa_0,J_0)$. Note that $g_0$ is simply the standard Euclidean
inner product.

Using the standard forms $\kappa_0$ and $\Omega_0$ by  straightforward computations
we can obtain some useful identities concerning $\kappa$-positive 3-forms.
\begin{lemma}
Let $(V,\kappa)$ be a symplectic vector space and $\Omega$ a
normalized $\kappa$-positive $3$-form, then we have
\begin{enumerate}
\item[1.] $\bigstar\Omega=-\Omega$ $($hence also $J\Omega=*\Omega)$;
\vspace{0.3cm}
\item[2.] $\Omega\wedge J\Omega=\frac{2}{3}\kappa^3$.
\end{enumerate}
\end{lemma}
\subsection{Decomposition of the exterior algebra}
Let $(V,\kappa)$ be an arbitrary 6-dimensional symplectic vector space and $\Omega$
a normalized $\kappa$-positive 3-form. Let us consider the natural
action of SU(3) on the exterior algebra
$\Lambda^{\bullet}V^*$.
Obviously SU(3) acts irreducibly on $V^*$ and $\Lambda^5V^*$, while
$\Lambda^2V^*$ and $\Lambda^3V^*$ decompose as follows:
\begin{equation}
\label{deco}
\begin{aligned}
&\Lambda^2V^*=\Lambda^2_1V^*\oplus\Lambda^2_6V^*\oplus\Lambda^2_8V^*\,,\\[3pt]
&\Lambda^3V^*=\Lambda^3_{Re}V^*\oplus
\Lambda^3_{Im}V^*\oplus\Lambda^{3}_6V^*\oplus
\Lambda^{3}_{12}V^*\,,
\end{aligned}
\end{equation}
where we set
\begin{itemize}
\item $\Lambda_1^2V^*=\R\,\kappa\,,$
\vspace{0.3 cm}
\item $\Lambda_6^2V^*=\{\bigstar(\alpha\wedge\Omega)\,|\,\alpha\in\Lambda^1V^*\}=
\{\f\in\Lambda^2V^*\,|\,J\f=-\f\}\,,$
\vspace{0.3 cm}
\item $\Lambda^2_8V^*=\{\f\in\Lambda^2V^*\,|\,\f\wedge\Omega=0 \mbox{ and
  }\bigstar{\f}=-\f\wedge\kappa\}$\\[3pt]
$\mbox{ }\mbox{ }\qquad \!=\{\f\in\Lambda^2V^*\,|\,J\f=\f\,,\,\f\wedge\kappa^2=0\}\,,\\$
\end{itemize}
and
\begin{itemize}
\item $\Lambda_{Re}^3 V^*=\R\, \Omega\,,$
\vspace{0.3 cm}
\item $\Lambda_{Im}^3 V^*=\R\,
J\Omega=\{\gamma\in\Lambda^3V^*\,|\,\gamma\wedge\kappa=0\,,\,\gamma\wedge\Omega=c\,\kappa^3,\,\,c\in\R\}\,,$
\vspace{0.3 cm}
\item
  $\Lambda_{6}^3V^*=\{\alpha\wedge\kappa\,|\,\alpha\in\Lambda^1V^*\}=\{\gamma
\in\Lambda^3V^*\,|\,\bigstar\gamma=\gamma\}\,,$
\vspace{0.3 cm}
\item
  $\Lambda^3_{12}V^*=\{\gamma\in\Lambda^3V^*\,|\,\gamma\wedge\kappa=0\,,\,\gamma\wedge\Omega=0\,,\,\gamma\wedge J\Omega=0\}\,.$
\end{itemize}
\vspace{0.5 cm}
\begin{rem}\label{emph}
\emph{
Now we emphasize some relations which will be useful:
\begin{enumerate}
\item[1.] If $\f\in\Lambda^{2}_6V^*\oplus\Lambda^{2}_{8}V^*$, then
  $\bigstar\f=-\f\wedge\kappa\,.$
\vspace{0.3 cm}
\item[2.] If
  $\gamma\in\Lambda^3_{Re}V^*\oplus\Lambda^3_{Im}V^*\oplus\Lambda^3_{12}V^*\,,$
  then $\bigstar \gamma=-\gamma$ and $\gamma\wedge\kappa=0$.
\vspace{0.3cm}
\item[3.] If $\alpha$ is an arbitrary 1-form, then
  $J(\alpha\wedge\Omega)=-\alpha\wedge\Omega$, consequently from the
  definition of $J$ it follows
$$
J\Omega\wedge\bigstar(\Omega\wedge\alpha)=-2\,\bigstar\alpha\,.
$$
\vspace{0.3cm}
\item[4.]
If $\beta \in \Lambda^2_8V^*$ then
\begin{eqnarray*}
*(\beta\wedge\beta)\wedge\kappa^2 & =
                  &\beta\wedge\beta\wedge *\kappa^2=2\,\beta\wedge\beta\wedge\kappa \\
                                  & = & -2\,\beta\wedge\bigstar\beta
                                    = -2|\beta|^2 \frac{\kappa^3}{6},
\end{eqnarray*}
so that
\begin{equation}
\label{ultima}
*(\kappa^2 \wedge *(\beta \wedge \beta))=-2|\beta|^2.
\end{equation}
\end{enumerate}
}
\end{rem}
We can obtain the decomposition of $\Lambda^4V^*$
using the duality given by the symplectic star operator.

Moreover we define the projections
$$
\begin{aligned}
&E_1\colon \Lambda^2V^*\to\Lambda^2_8V^*\,,\\
&E_2\colon \Lambda^3V^*\to\Lambda^3_{12}V^*
\end{aligned}
$$
by
\begin{eqnarray}
\label{E_1}
&&E_1(\alpha)=\frac{1}{2}(\alpha+J\alpha)-
\frac{1}{18}*((*(\alpha+J\alpha)+(\alpha+J\alpha)\wedge \kappa)\wedge\kappa)\,\kappa\,,\\
\label{E_2}
&&E_2(\beta)=\beta-\frac{1}{2}*(J\beta\wedge\kappa)\wedge\,\kappa
               -\frac{1}{4}*(\beta\wedge J\Omega)\,\Omega
               -\frac{1}{4}*(\Omega\wedge\beta) \,J\Omega\,.
\end{eqnarray}

Note that $E_2$ commutes with $*$ since the latter is an automorphism of 
$\Lambda^3_{12}V^*$. The same is true for $J$ (hence also for $\bigstar$).
\subsection{The $\e$-identities}
As done by Bryant in the
G$_2$ case we introduce the following $\e$-notation, which will be
useful in the sequel.
$$
\Om_0=\frac{1}{6}\e_{ijk}\,e^{ijk}\,,\quad*\Om_0=\frac{1}{6}\ov{\e}_{ijk}\,e^{ijk}\,,\quad
\kappa_0=\frac{1}{2}\kappa_{ij}\,e^{ij}\,.
$$
We  will use
the  following identities, whose proof is straightforward:
\begin{equation}
\label{ep}
\begin{aligned}
&\e_{ipq}\kappa_{pq}=0\,;\\[3pt]
&\kappa_{ip}\kappa_{pj}=-\delta_{ij}\,;\\[3pt]
&\e_{ijp}\kappa_{pr}=\ov{\e}_{ijr}\,;\\[3pt]
&\ov{\e}_{ijp}\kappa_{pr}=-\e_{ijr}\,;\\[3pt]
&\ov{\e}_{ipq}\e_{jpq}=-4\kappa_{ij}\,;\\[3pt]
&\e_{ipq}\e_{jpq}=4\delta_{ij}=\ov{\e}_{ipq}\ov{\e}_{jpq}\,;\\[3pt]
&\ov{\e}_{ijp}\e_{klp}=-\kappa_{ik}\delta_{jl}+\kappa_{jk}\delta_{il}+
\kappa_{il}
\delta_{jk}-\kappa_{jl}\delta_{ik}\,;\\
&\e_{ijp}\e_{klp}=-\kappa_{ik}\kappa_{jl}+\kappa_{il}\kappa_{jk}+\delta_{ik}
\delta_{jl}-\delta_{jk}\delta_{il}=\ov{\e}_{ijk}\ov{\e}_{ipq}\,.\\
\end{aligned}
\end{equation}
These equations will be called $\epsilon$-\emph{identities}.
As a  first application of these formulae we can decompose the Lie
algebra $\mathfrak{so}(6)$ as follows. Consider the real
representation of complex matrices induced by $J_{0}$
$$
\rho\colon \mathfrak{gl}(3,\C)\to\mathfrak{gl}(6,\R)\,,
$$
where $\rho(A)$ is the block matrix $(B_{ij})_{i,j=1,2,3}$, with
$
B_{ij}=\left(
\begin{array}{cc}
\mbox{Re}\, a_{ij} & \mbox{Im}\, a_{ij}\\
-\mbox{Im}\, a_{ij} & \mbox{Re}\, a_{ij}
\end{array}
\right).
$
Thus a matrix $A=(a_{ij})$ lies in $\mathfrak{su}(3)$ if and only if
$$
\e_{ijk}a_{jk}=0\quad\mbox{and}\quad \kappa_{jk}a_{jk}=0\,.
$$
So we have the decomposition
\begin{equation*}
\mathfrak{so}(6)=\mathfrak{su}(3)\oplus[\mathbb{R}]_1\oplus[\mathbb{R}^6]_2,
\end{equation*}
where
$$
([a]_1)_{ij}=a\,\kappa_{ij}\,,\qquad ([v]_2)_{ij}=\e_{ijp}\,v_p\,.
$$
\subsection{Decomposition of symmetric 2-tensors}
\label{symmetric}
In order to express the Ricci tensor in terms of skew-symmetric
forms we must establish the correspondence which we are going to
describe. The 21-dimensional space of symmetric covariant 2-tensor
on $V$ splits into irreducible $\mathfrak{su}(3)$-modules as
follows:
$$
S^2V^*=\R\,g_0\oplus S_+^2\oplus S_-^2\,,
$$
where
$$
\begin{aligned}
&S_+^2=\{h\in S^2V^*\,:\,J_0h=h,\,\mbox{tr}_{g_0}h=0 \}\,,\\
&S_-^2=\{h\in S^2V^*\,:\,J_0h=-h\}\,.
\end{aligned}
$$
We will denote by $S_0^2$ the direct sum $S_+^2\oplus S_-^2$.

The maps
$$
\begin{aligned}
&\iota\colon S^2_+\longrightarrow \Lambda_{8}^2V^*\,,\\
&\gamma\colon S^2_-\longrightarrow \Lambda_{12}^3V^*
\end{aligned}
$$
defined by
$$
\begin{aligned}
\iota(h_{ij}e^ie^j)&=h_{ip}\kappa_{pj}\,e^{ij}\,,\\
\gamma(h_{ij}e^ie^j)&=h_{ip}\epsilon_{pjk}\,e^{ijk}\,.
\end{aligned}
$$
are isomorphisms of $\mathfrak{su}(3)$-representations.
\subsection{SU(3)-structures on manifolds.}
\label{SU3man}
Let $M$ be a $6$-dimensional manifold.
A SU(3)-structure on $M$ is determined by  the choice of:
\begin{itemize}
\item a non-degenerate 2-form $\kappa$,
\vspace{0.1 cm}
\item a normalized $\kappa$-positive 3-form $\Omega$ (i.e. $\Omega[x]$ is
  $\kappa[x]$-positive and normalized at every $x$ in $M$).
\end{itemize}
In fact, as we have seen, $\Omega$ determines a $\kappa$-calibrated almost complex
structure $J$ such that $\varepsilon=\Omega+iJ\Omega$ is of type
$(3,0)$ and satisfies  equation \eqref{cienne}.
We refer to $\varepsilon$ as to the \emph{complex volume of $(\kappa,\Omega)$}.
In the sequel the induced scalar product will be denoted by $g$ or
alternatively by $\langle\,,\,\rangle$ and the associated Hodge
operator by $*$.\\
Note that the SU$(3)$-structure determined by
$(\kappa,\Omega)$ is integrable if and only if
\begin{equation}
\label{integrabili}
d\kappa=0\,,\quad d\Omega=d^*\Omega=0\,.
\end{equation}
In fact, since $J\Omega=*\Omega$, equations \eqref{integrabili}
are equivalent to
$$
d\kappa=0\,,\quad d\varepsilon=0\,.
$$
Hence, since $\varepsilon\wedge\ov{\varepsilon}=i\frac{4}{3}\,\kappa^{3}$,
remark \ref{rem1} implies
$$
\nabla\kappa=0\,,\quad \nabla J=0\,,\quad \nabla\varepsilon=0\iff
d\kappa=0\,,\quad d\varepsilon=0\,.\quad\mbox{          }
$$
\subsection{Torsion forms}
\label{torsionforms}
Let $(M,\kappa,\Omega)$ be a SU$(3)$-manifold.
According with \eqref{deco} the spaces of $r$-forms splits in
$\mathfrak{su}(3)$-modules as follows:
\begin{equation*}
\begin{aligned}
&\Lambda^2M=\Lambda^2_1M\oplus\Lambda^2_6M\oplus\Lambda^2_8M\,,\\[3 pt]
&\Lambda^3M=\Lambda^3_{Re}M\oplus
\Lambda^3_{Im}M\oplus\Lambda^{3}_6M\oplus
\Lambda^{3}_{12}M\,,\\[3 pt]
&\Lambda^4M=\Lambda^4_1M\oplus\Lambda^4_6M\oplus\Lambda^4_8M\,,\\
\end{aligned}
\end{equation*}
where the meaning of symbols is obvious. Consequently the derivatives
of the structure forms decompose as
\begin{equation}
\begin{aligned}
\label{scomposition}
d\kappa&=\nu_0\,\Omega+\alpha_0\,J\Omega+\nu_1\wedge\kappa+\nu_3\,,\\
d\Omega& =\pi_0\,\kappa^2+\pi_1\wedge\Omega-\pi_2\wedge\kappa\,,\\
dJ\Omega&=\sigma_0\,\kappa^2+\sigma_1\wedge\Omega-\sigma_2\wedge\kappa\,,
\end{aligned}
\end{equation}
where $\nu_0,\alpha_0,\pi_0,\sigma_0\in C^{\infty}(M,\R)$,
$\nu_1,\pi_1,\sigma_1\in \Lambda^1M$, $\pi_2,\sigma_2\in\Lambda^2_8M$
and $\nu_3\in\Lambda^3_{12}Mù$.
\newline

\noindent The following equations are derived from a G$_2$ formula which was
obtained in \cite{B1}.
\begin{lemma}
\label{Bryantlemma}
With the notations introduced above
\begin{eqnarray}
\label{Bryant 2}
&& J\Omega\wedge(*dJ\Omega)-(*d\Omega)\wedge\Omega=0\,.
\end{eqnarray}
\end{lemma}
\begin{proof}
See the appendix.
\end{proof}
\noindent Now we are able to prove the following
\begin{theorem}
\label{relazionetorsioni}
The following relations hold:
\begin{enumerate}
\item[1.]$\pi_0=\frac{2}{3} \alpha_0$\,,
\item[2.]$\sigma_0=-\frac{2}{3} \nu_0$\,,
\item[3.]$\sigma_1=J\pi_1$\,.
\end{enumerate}
\end{theorem}
\begin{proof}
1. From the relation $\Omega\wedge\kappa=0$ it follows
$$
\begin{aligned}
0&=d(\Omega\wedge\kappa)=d\Omega\wedge\kappa-\Omega\wedge d\kappa \\
&=\pi_0\,\kappa^3-\pi_2\wedge\kappa^2-\alpha_0\,\Omega\wedge J\Omega-\Omega\wedge\nu_3\\
&=(\pi_0-\frac{2}{3}\alpha_0)\,\kappa^3\,,
\end{aligned}
$$
where we have used that $\pi_2\wedge\kappa^2=0$, $\Omega\wedge\nu_3=0$.\\
2. Analogous to 1 starting from $\kappa\wedge J\Omega=0$.\\
3. This formula is a consequence of formula \eqref{Bryant 2} together
  with the definition of $J$. We have
$$
\begin{aligned}
0&=(*d\Omega)\wedge\Omega-J\Omega\wedge*dJ\Omega\\
 &=*(\pi_1\wedge\Omega)\wedge\Omega-J\Omega\wedge*(\sigma_1\wedge\Omega)\\
 &=-J(\bigstar(\pi_1\wedge \Omega)\wedge J\Omega)-J(\Omega\wedge\bigstar
 (\sigma_1\wedge\Omega))\\
&=J(J\Omega\wedge\bigstar(\Omega\wedge\pi_1))+J(\Omega\wedge\bigstar
 (\Omega\wedge\sigma_1))\,.
\end{aligned}
$$
Applying the definition of $J$ and remark \ref{emph} we get
$$
J(-2\bigstar\pi_1)-J(2J\bigstar\sigma_1)=-2J\bigstar\pi_1+2\bigstar\sigma_1=0\,,
$$
i.e.
$$
\sigma_1=J\pi_1\,.
$$
\end{proof}
Hence we can rewrite \eqref{scomposition} as:
$$
\begin{aligned}
d\kappa&=-\frac{3}{2}\sigma_0\,\Omega+\frac{3}{2}\pi_0\,J\Omega+\nu_1\wedge\kappa+\nu_3\,;\\[3pt]
d\Omega& =\pi_0\,\kappa^2+\pi_1\wedge\Omega-\pi_2\wedge\kappa\,;\\[3pt]
dJ\Omega&=\sigma_0\,\kappa^2+J\pi_1\wedge\Omega-\sigma_2\wedge\kappa\,.
\end{aligned}
$$
\begin{definition}
The forms $\{\pi_0,\sigma_0,\pi_1,\nu_1,\sigma_2,\nu_3\}$ are called the
\emph{torsion forms} of the \emph{SU(3)}-structure.
\end{definition}
A SU(3)-structure is
integrable if and only if all of the torsion forms vanish
identically.
\newline

Several interesting
special SU$(3)$-structures can be described in terms of torsion forms.
\begin{enumerate}
\item[1.]\textbf{6-dimensional GCY structures.} let $(M,\kappa,\Omega)$ be a
  6-dimensional GCY manifold. The equation $d\kappa=0$ implies
$$
\pi_0=\sigma_0=0\,,\quad \nu_1=0\,,\quad\nu_3=0\,.
$$
Therefore $d\Omega$ and $dJ\Omega$  reduce to
$$
\begin{aligned}
d\Omega&=\pi_1\wedge\Omega-\pi_2\wedge\kappa\,,\\
dJ\Omega&=J\pi_1\wedge\Omega-\sigma_2\wedge\kappa\,.
\end{aligned}
$$
Since the complex volume form $\varepsilon$ associated to $(\kappa,\Omega)$
is of type (3,0), $\overline{\partial}_J \varepsilon$ is the $(3,1)$-part
(hence the $J$ anti-invariant part) of $d\varepsilon$. Thus we have
$$
\overline{\partial}_J\varepsilon=\frac{1}{2}(d\varepsilon-Jd\varepsilon)\,.
$$
Thus
$$
\begin{aligned}
\overline{\partial}_J\varepsilon=&\frac{1}{2}(d\varepsilon-Jd\varepsilon)\\
=&\frac{1}{2}(d\Omega+idJ\Omega-Jd\Omega-iJdJ\Omega)\\
=&\frac{1}{2}\{d\Omega-Jd\Omega+i(dJ\Omega-JdJ\Omega)\}\\
=&\frac{1}{2}\{\pi_1\wedge\Omega-J(\pi_1\wedge\Omega)
+i(J\pi_1\wedge\Omega-J(J\pi_1\wedge\Omega))\}\\
=&\pi_1\wedge\Omega+i\,J\pi_1\wedge\Omega.
\end{aligned}
$$
Hence  by lemma \ref{facile} the equation $\overline{\partial}_J\varepsilon=0$ is equivalent
to $\pi_1=0$. It follows that 6-dimensional GCY structures can be defined as
SU(3)-structures satisfying
$$
\pi_0=\sigma_0=0\,,\quad \nu_1=\pi_1=0\,,\quad\nu_3=0\,.
$$
\item[2.]\textbf{Special generalized Calabi-Yau structure.} These
  structures has been introduced and studied first by P. de
  Bartolomeis in \cite{dB}.
\begin{definition} Let $M$ be a $6$-dimensional manifold.
A special generalized Calabi-Yau
structure \emph{(SGCY)} on $M$ is a \emph{SU(3)}-structure such that
the defining forms
$\kappa$, $\Omega$ are closed, i.e.
$$
d\kappa=0\,,\quad
d\Omega=0\,.
$$
\end{definition}
Special  generalized Calabi-Yau manifolds can be considered as a
subclass of generalized Calabi-Yau manifold, in fact it is
immediately verified that in this case the complex volume form
$\varepsilon$ associated to $(\kappa,\Omega)$ satisfies the
condition 2 of definition \ref{GCY} (see \cite{dBT}). SGCY
manifolds are taken into consideration also in \cite{Bed}, \cite{CT}
and \cite{TV}.\\
Such a structure can be characterized by
$$
\pi_0=\sigma_0=0\,,\quad\nu_1=\pi_1=0\,,\quad\pi_2=0\,,\quad\nu_3=0\,.
$$
\item[3.]\textbf{Half-flat structure.} Half-flat manifolds have a
central role in the evolution theory
developed by Hitchin in \cite{Hi2}  and can be used to construct
non-compact examples of G$_2$-manifolds.
\begin{definition}
A \emph{SU}$(3)$-structure $(\kappa,\Omega)$ is said to be \emph{half-flat}
if the structure forms satisfy the equations
$$
d(\kappa\wedge\kappa)=0\,,\quad d\Omega=0\,.
$$
\end{definition}
Let $(\kappa,\Omega)$ be a half-
flat structure. By the hypothesis $d\Omega=0$ we get
$$
\pi_{i}=0\,,\quad i=0,1,2\,;
$$
then
$$
d\kappa=-\frac{3}{2}\,\sigma_0\,\Omega+\nu_1\wedge\kappa+\nu_3\,.
$$
On the other hand the hypothesis $d(\kappa\wedge\kappa)=0$ implies
$$
0=d\kappa\wedge\kappa=-\frac{3}{2}\sigma_0\,\Omega\wedge\kappa+\nu_1\wedge\kappa^2+\nu_3\wedge\kappa=\nu_1\wedge\kappa
^2\,,
$$
which forces $\nu_1$ to vanish, since the exterior multiplication by
$\kappa^2$ is an isomorphism on $\Lambda^1M$.
Therefore half-flat structures can be described as SU(3)-structures
satisfying
$$
\pi_{i}=0\,,\quad i=0,1,2\,,\quad \nu_1=0\,.
$$
\end{enumerate}
\subsection{Some SU$(3)$ representation theory}
Every irreducible representation $\rho$ of the simple Lie group
SU$(3)$ can be labeled by a pair of integers $(p,q)$ that represent
the highest weight of $\rho$ with respect to a fixed base of the
root system of a fixed maximal torus of SU$(3)$. We will denote
$\rho$ by $\lambda_{p,q}$. Nevertheless in the sequel we need to
deal with {\it
  real}
representation of SU$(3)$, so (similar as in \cite{S}) we will define the
irreducible real representations $V_{p,q}$ ($p\neq q$) and $V_{p,p}$ by
$$
\begin{aligned}
&V_{p,q}\otimes_{\R}\C=\lambda_{p,q}\oplus\lambda_{q,p}\,,\\
&V_{p,p}\otimes_{\R}\C=\lambda_{p,p}\,.
\end{aligned}
$$
Keeping in mind this fact, we can use the complex representation theory
to decompose a given real SU$(3)$-representation into irreducible real
SU$(3)$-modules. As it is well-known (see \cite{B2}) the polynomial pointwise
invariants of order $k$ are polynomials in a canonically defined
section of the vector bundle
$$
\mathcal{Q}\times_{\rho_{1}\times\dots\times\rho_{k}}(V_1(\mathfrak{su}(3))\oplus\dots\oplus V_{k}(\mathfrak{su}(3)))\,,
$$
where $\mathcal{Q}$ is the SU$(3)$-reduction and
$V_{j}(\mathfrak{su}(3))$ is the SU$(3)$-representation uniquely
defined by
$$
(\mathfrak{gl}(6,\R)/\mathfrak{su}(3))\otimes
S^{j}(\R^6)=V_{j}(\mathfrak{su}(3))\oplus(\R^6\otimes S^{j+1}(\R^6))\,.
$$
For the first order invariants we have
$$
V_{1}(\mathfrak{su}(3))=\mathfrak{so}(6)/\mathfrak{su}(3)\otimes \R^6
$$
so that
$$
V_1(\mathfrak{su}(3))=2\,V_{0,0}\oplus 2\,(\R^6)^*\oplus 2\,\Lambda^2_8\oplus
\Lambda^3_{12}
$$
which matches with to the degree and types of our torsion forms.
Rather standard calculation in $\mathfrak{su}(3)$-representation
theory allow us to decompose also the 252-dimensional representation
$V_{2}(\mathfrak{su}(3))$ into $\mathfrak{su}(3)$-irreducible
submodules
$$
\begin{aligned}
V_{2}(\mathfrak{su}(3))=&3\,V_{0,0}\oplus 4\,V_{1,0}\oplus
5\,V_{1,1}\oplus 3\,V_{2,1}\oplus 4\,V_{2,0} \oplus V_{3,0}\oplus V_{2,2ù}
\,,
\end{aligned}
$$
\section{Riemannian invariants of SU(3)-structures}
\label{Riemanniansection}
\subsection{The Levi-Civita connection}
Fix a SU(3)-reduction $\mathcal{Q}$ of the linear frame bundle $\mathcal{L}(M)$, given
by the pair $(\kappa,\Omega)$.  $\mathcal{Q}$ is a subbundle of the principal
SO(6)-bundle $p\colon\mathcal{F}\to M$ of the normal frames of the metric $g$ associated
to the pair $(\kappa,\Omega)$. Consider on the bundle $\mathcal{F}$
the tautological $\R^6$-valued $1$-form $\omega$ defined by
$\omega[u](v)=u(p_{*}[u]v)$ for every $u\in \mathcal{F}$ and $v\in T_u\mathcal{F}$.
On $\mathcal{F}$ we have also the Levi-Civita connection 1-form
$\psi$ taking values in $\mathfrak{so}(6)$. Using the canonical basis
$\{e_1,\dots,e_6\}$ of $\R^{6}$ we will regard
$\omega$ as a vector of $\R$-valued 1-forms on $\mathcal{F}$
$$
\omega=\omega_1e_1+\dots+\omega_6e_6
$$
and $\psi$ as a skew-symmetric
matrix of 1-forms, i.e. $\psi=(\psi_{ij})$. With these notations the
first structure equation relating $\omega$ and $\psi$
\begin{equation}
\label{1se}
d\omega=-\psi\wedge\omega\,,
\end{equation}
becomes $d\omega_i=-\psi_{ij}\wedge\omega_j\,$. Note that equation  (\ref{1se})
simply means that $\psi$ is torsion-free.\\
The curvature of $\psi$ is by definition the $\mathfrak{so}(6)$-valued
2-form $\Psi=d\psi+\psi\wedge\psi$. In index notation
$$
\Psi_{ij}=d\psi_{ij}+\psi_{ik}\wedge\psi_{kj}=\frac{1}{2}R_{ijkl}\,\omega_k\wedge\omega_{l}\,.
$$

We consider the pull-backs of $\psi$ and $\omega$ to
$\mathcal{Q}$  and denote them by the same symbols for the sake of brevity.
The intrinsic torsion of the SU(3)-structure
measures the failing of $\psi$ to take values in $\mathfrak{su}$(3). More
precisely, according to the splitting $\mathfrak{so}(6)=\mathfrak{su}(3)\oplus[\mathbb{R}]_1\oplus
[\mathbb{R}^6]_2$, we decompose $\psi$ as follows
\begin{equation*}
\psi=\theta+[\mu]_1+[\tau]_2\,.
\end{equation*}
Thus $\theta$ is a connection 1-form on $\mathcal{Q}$ which in general
is not torsion-free.\\
As before we shall regard $\tau$ as a vector of 1-forms
$\tau=\tau_ie_i$. Furthermore we can write
\begin{equation}
\label{deftaumu}
\tau_{i}=T_{ij}\,\omega_j\quad\mbox{and}\quad \mu=M_{i}\,\omega_i\,,
\end{equation}
where $T_{ij}$ and $M_i$ are smooth functions. The fact that $\psi$
is torsion-free implies
\begin{equation}
\label{domega}
d\omega_i=-\theta_{ij}\wedge\omega_j-\epsilon_{ijk}\,\tau_{k}\wedge\omega_j-\kappa_{ij}\,\mu\wedge\omega_j
\,.
\end{equation}
\subsection{The curvature in index notation}
In order to decompose the curvature 2-form we give the following
\begin{lemma}
\label{lemmaluigi}
These identities hold:
\begin{enumerate}
\item[1.]$\theta\wedge[\mu]_1+[\mu]_1\wedge\theta=0$ ;
\vspace{0.1 cm}
\item[2.]$[\tau]_2\wedge[\mu]_1-[\mu]_1\wedge[\tau]_2=0$ ;
\vspace{0.1 cm}
\item[3.]$\theta\wedge[\tau]_2+[\tau]_2\wedge\theta=[\theta\wedge\tau]_2$ ;
\vspace{0.1 cm}
\item[4.]$[\tau]_2\wedge[\mu]_1+[[\mu]_1\wedge\tau]_2=0$ .
\end{enumerate}
\end{lemma}
\begin{proof}
The proof is a straightforward application of $\epsilon$-identities \eqref{ep}.
To see how things work, we prove the first one. Since $\theta$ takes
values in $\mathfrak{su}(3)$ we have
\begin{equation*}
\epsilon_{pkl}\,\theta_{kl}=\epsilon_{klp}\,\theta_{kl}=0\,.
\end{equation*}
So
\begin{equation*}
\overline{\epsilon}_{ijp}\epsilon_{klp}\,\theta_{kl}=0
\end{equation*}
for every $i,j=1,\dots,6$. Then applying the $\epsilon$-identities
\eqref{ep} we get
\begin{equation*}
\begin{aligned}
0=&\overline{\epsilon}_{ijp}\epsilon_{klp}\,\theta_{kl}\\
 =&(-\kappa_{ik}\delta_{jl}+\kappa_{jk}\delta_{il}+
\kappa_{il}
\delta_{jk}-\kappa_{jl}\delta_{ik})\,\theta_{kl}\\
=&2\kappa_{jk}\,\theta_{ki}-2\kappa_{ik}\,\theta_{kj}\,,
\end{aligned}
\end{equation*}
i.e.
\begin{equation*}
\kappa_{jk}\,\theta_{ki}=\kappa_{ik}\,\theta_{kj}\,.
\end{equation*}
Consequently
$$
\theta_{ik}\wedge \kappa_{kj}\,\mu+\kappa_{ik}\,\mu\wedge\theta_{kj}=0\,,
$$
i.e.
$$
\theta\wedge[\mu]_1+[\mu]_1\wedge\theta=0\,.
$$
\end{proof}
Now we can introduce the following quantities
\begin{eqnarray}
&&D\theta=d\theta+\,\theta\wedge\theta+\,[\tau]_2\wedge[\tau]_2-\frac{2}{3}[\kappa_{ij}\,\tau_i
\wedge\tau_j]_1\,,\\
\label{Dtau}
&& D\tau=d\tau+\,\theta
\wedge\tau-2\,[\mu]_1\wedge \tau\,,\\
\label{Dmu}
&&D\mu=d\mu+\frac{2}{3}\kappa_{ij}\,\tau_i\wedge\tau_j\,.
\end{eqnarray}
With this definition $D\theta$ takes values in
$\mathfrak{su}(3)$. Moreover by lemma \ref{lemmaluigi} we get
\begin{equation*}
\begin{aligned}
\Psi=& d(\theta+[\tau]_2+[\mu]_1)+(\theta+[\tau]_2+[\mu]_1)\wedge
(\theta+[\tau]_2+[\mu]_1)\\
=&D\theta+[D\tau]_2+[D\mu]_1\,.
\end{aligned}
\end{equation*}
Using the $\omega$-frame we shall write
\begin{eqnarray}
&&D\theta_{ij}=\frac{1}{2}S_{ijkl}\,\omega_k\wedge\omega_l\,,\\
\label{full}
&&D\tau_i =\frac{1}{2}T_{ijk}\,\omega_j\wedge\omega_k\,,\\
\label{coppia}
&&D\mu=\frac{1}{2}N_{kl}\,\omega_k\wedge\omega_l\,.
\end{eqnarray}
By the definition of the curvature form we have
\begin{equation*}
R_{ijkl}=S_{ijkl}+\epsilon_{ijp}T_{pkl}+\kappa_{ij}N_{kl}\,.
\end{equation*}
In this notation the first Bianchi identity
$$
\Psi\wedge\omega=0\,,
$$
has the indicial expression
\begin{equation}
\begin{aligned}
\label{Bianchi}
&S_{ijkl}+S_{iljk}+S_{iklj}+\\
&+\epsilon_{ijp}T_{pkl}+\epsilon_{ilp}T_{pjk}+\epsilon_{ikp}T_{plj}
+\kappa_{ij}N_{kl}+\kappa_{il}N_{jk}+\kappa_{ik}N_{lj}=0
\end{aligned}
\end{equation}
Let $Ric_{ij}=R_{ikkj}$ and $s=Ric_{kk}$  be respectively the Ricci
tensor and the scalar curvature of $(M,g)$. Starting from equation
\eqref{Bianchi} a long, but straightforward computation gives the
following
\begin{theorem}
\label{controbuio}
In the previous notation we have
\begin{eqnarray*}
&&Ric_{ij}=2\epsilon_{ipq}T_{pqj}-3\kappa_{ip}N_{pj}\,,\\[3pt]
&&s=2\epsilon_{kpq}T_{pqk}-3\kappa_{kp}N_{pk}\,.
\end{eqnarray*}
\end{theorem}
\subsection{Ricci tensor in terms of torsion forms} Denote by $\pi$ the projection
$\pi\colon\mathcal{Q}\to M$. In terms of the $\omega$-frame the pull-backs of
the structure forms take their standard expression, i.e.
\begin{eqnarray*}
&&\pi^*(\Omega)=\frac{1}{6}\epsilon_{ijk}\,\omega_i\wedge\omega_j\wedge\omega_k\,,\\
&&\pi^*(J\Omega)=\frac{1}{6}\overline{\epsilon}_{ijk}\,\omega_i\wedge\omega_j\wedge\omega_k\,,\\
&&\pi^*(\kappa)=\frac{1}{2}\kappa_{ij}\,\omega_i\wedge\omega_j\,.
\end{eqnarray*}
Taking into account formula \eqref{domega} and $\epsilon$-identities,
we immediately get
\begin{prop}
The derivatives of the structure forms are
\begin{equation*}
\begin{aligned}
&d\pi^*(\Omega)=\frac{1}{2}(-\kappa_{ja}\kappa_{kb}+\kappa_{jb}\kappa_{ka})\,\tau_{b}
\wedge\omega_{a}\wedge\omega_j\wedge
\omega_k-3\,\mu\wedge \pi^*(J\Omega)\,,\\[3pt]
&d\pi^*(J\Omega)=(\tau_{j}\wedge\omega_j)\wedge\pi^*(\kappa)-3\,\mu\wedge\pi^*(\Omega)\,,\\[3
pt]
&d\pi^*(\kappa)=\overline{\epsilon}_{lrj}\,\tau_{l}\wedge\omega_{r}\wedge\omega_{j}\,.
\end{aligned}
\end{equation*}
\end{prop}
Now we can decompose the derivatives of the structure forms: a direct
computation gives the following formulae
\begin{equation*}
\begin{aligned}
\pi^*(\pi_0)&=\frac{2}{3}T_{ii}\,,\\
\pi^*(\pi_1)&=\e_{ijk}T_{ij}\,\omega_k+3\kappa_{ik}M_{i}\,\omega_k\,,\\
\pi^*(\pi_2)&=\frac{1}{2}\ov{\e}_{sra}\e_{aij}T_{sr}\,\omega_{i}\wedge\omega_j
-2\kappa_{ia}T_{aj}\,\omega_{i}\wedge\omega_j+\frac{2}{3}T_{ii}\,\pi^*(\kappa)\,,
\\
\pi^*(\sigma_0)&=\frac{2}{3}\kappa_{ij}\,T_{ij}\,,\\
\pi^*(\sigma_2)&=\frac{1}{2}\e_{rsa}\e_{aij}T_{rs}\,\omega_i\wedge\omega_j-2T_{ij}
\omega_i\wedge\omega_j+\frac{2}{3}\kappa_{ij}T_{ij}\,\pi^*(\kappa)\,,\\
\pi^*(\nu_1)&=\e_{ijk}T_{ij}\,\omega_k\,,\\
\pi^*(\nu_3)&=\ov{\e}_{aij}T_{ak}\,\omega_i\wedge\omega_j\wedge\omega_k
+\frac{1}{6}\kappa_{ab}T_{ab}\e_{ijk}\,\omega_i\wedge\omega_j\wedge\omega_k\\
&\quad-\frac{1}{6}T_{aa}\ov{\e}_{ijk}\,\omega_i\wedge\omega_j\wedge\omega_k
-\frac{1}{2}T_{ab}
\e_{abi}\kappa_{jk}\,\omega_i\wedge\omega_j\wedge\omega_k\,.
\end{aligned}
\end{equation*}
\textbf{Warning:}
From now on we identify the torsion forms with their pull-backs to the
principal SU$(3)$-bundle $\mathcal{Q}$.\\
\newline
Combining the previous formulae and \eqref{domega} we are able to prove the following (see
the appendix)
\begin{theorem}
\label{t s}
In terms of torsion forms  the scalar curvature of the metric induced
by the \emph{SU(3)}-structure is
expressed  as
\begin{equation}
\label{curvatura scalare}
\begin{aligned}
s=&\frac{15}{2}\pi_0^2+\frac{15}{2}\sigma_0^2+2d^*\pi_1+2d^*\nu_1
-|\nu_1|^2-\frac{1}{2}|\sigma_2|^2\\
&-\frac{1}{2}|\pi_2|^2
-\frac{1}{2}|\nu_3|^2+4\langle\pi_1,\nu_1\rangle\,.
\end{aligned}
\end{equation}
\end{theorem}
Here we collect some consequences of formula (\ref{curvatura scalare})
when the SU(3)-structure has special features.
\begin{enumerate}
\item[1.]\textbf{GCY structure.} The condition
  $\overline{\partial}_J\epsilon=0$ reads as $\pi_1=0$ (see section
  \ref{torsionforms}), so that, taking into account $d\kappa=0$,
$$
s=-\frac{1}{2}|\sigma_2|^2-\frac{1}{2}|\pi_2|^2\,.
$$
\item[2.]\textbf{SGCY structure.} This is a special case of the previous
one with the extra-condition $\pi_2=0$. The scalar curvature takes the form
\begin{equation}
\label{curvatura scalare SGCY}
s=-\frac{1}{2}|\sigma_2|^2\,.
\end{equation}
\item[3.] \textbf{Half-flat structure.} The condition $d\kappa\wedge
  \kappa=0$ reads in terms of torsion forms as $\nu_{1}=0$. Thus in the
  half-flat case the
  scalar curvature takes the form
$$
s=\frac{15}{2}\sigma_0^2-\frac{1}{2}|\sigma_2|^2-\frac{1}{2}|\nu_3|^2\,.
$$
\end{enumerate}
\begin{cor}
The scalar curvature of a $6$-dimensional generalized Calabi-Yau manifold is
everywhere non-positive and it vanishes identically if and only if the
\emph{SU(3)}-structure has no torsion.
\end{cor}
Now we write the Ricci curvature
$Ric_{ij}=2\epsilon_{ipq}T_{pqj}-3\kappa_{ip}N_{pj}$ in terms of the
torsion forms using the operators $\iota$ and $\gamma$ defined in
section \ref{symmetric}.
\begin{theorem}
\label{teoremaricci}
If $M$ is endowed with the \emph{SU}$(3)$-structure $(\kappa,\Omega)$ with
torsion forms given by \eqref{scomposition} , then the traceless part of the Ricci
tensor of the induced metric is
\begin{equation}
\label{riccitensor}
Ric_0=\iota^{-1}(E_1(\phi_1))+\gamma^{-1}(E_2(\phi_2))\,,
\end{equation}
where
\begin{equation*}
\begin{aligned}
\phi_1=&-*(\nu_1\wedge J\nu_3)+\frac{1}{4}\,*(\pi_2\wedge\pi_2)+
        \frac{1}{4}\,*(\sigma_2\wedge\sigma_2)+\\
       &+dJ\pi_1+\frac{1}{2}\,d^*\nu_3
        +\frac{1}{2}\,d^*(\nu_1\wedge\kappa)-\frac{1}{4}\,d*(\pi_0\,\Omega)
        +\frac{1}{4}\,d^*(\sigma_0\,\Omega)\,,\\[10pt]
\phi_2=&-2\sigma_0\,\nu_3-4\,\sigma_2\wedge\nu_1-2\,Jd\pi_2-2\,\bigstar d\sigma_2
        -4\,d*(\nu_1\wedge*\Omega)+\\[4pt]
       &-2d*(J\pi_1\wedge\Omega)+2\pi_0\,J\,\nu_3-2\,Jd*(\pi_1\wedge\Omega)-
        4\,\pi_2\wedge J\pi_1+\\
       &+4\,\nu_1\wedge*(J\pi_1\wedge\Omega)-2\,J\nu_1\wedge*(\nu_1\wedge\Omega)
        -\frac{1}{2}Q(\nu_3,\nu_3)\,,
\end{aligned}
\end{equation*}
$E_1$ and $E_2$ are the maps defined by equations \eqref{E_1} and
\eqref{E_2}
and $Q$ is the bilinear form
$Q\colon\Lambda^3_{12}M\times\Lambda^3_{12}M\to\Lambda^3M$
defined by
$$
Q(\alpha,\beta)=\e_{ijl}\iota_{e_j}\iota_{e_i}\,\alpha\wedge\iota_{e_l}\beta\,,
$$
where $\{e_1,\dots,e_6\}$ is a unitary frame and $\iota$ denotes the
contraction of forms.
\end{theorem}
\begin{rem}\emph{
The formulae for the scalar curvature and for the traceless part of
the Ricci tensor are justified by representation theory. Both $s$ and
$Ric_0$ must be the linear combination of linear terms in
$V_2(\mathfrak{su}(3))$ and quadratic terms in
$V_1(\mathfrak{su}(3))$. For the scalar curvature the terms must take values
in the $V_{0,0}$ copies of $V_1$ and $V_2$, while for the Ricci curvature
the terms must take values  in $\Lambda^2_8$ and $\Lambda^3_{12}$ copies
of $V_1$ and $V_2$. (For $S^2_0=\Lambda^2_8\oplus \Lambda^3_{12}$). So
we have to consider:
$$
\begin{aligned}
S^2(V_1(\mathfrak{su}(3)))=&11\,V_{0,0}\oplus 13\,V_{1,0}
\oplus 17\, V_{1,1}\oplus 12\,V_{2,0}\oplus\\
&\oplus 3\,V_{3,0}\oplus 4\,V_{2,2}\oplus 9\,
V_{2,1}\oplus 2\,V_{3,1}\,.
\end{aligned}
$$
The 11 copies of $V_{0,0}$ are generated by
\begin{itemize}
\item $\pi_0^2,\,\sigma_0^2,\,\pi_0\sigma_0$;
\vspace{0,2cm}
\item $\,|\pi_1|^2,\,|\nu_1|^2,\,<\pi_1,\nu_1>$ and another bilinear
expression in $\pi_1$, $\nu_1$ which does not
appear in formula \eqref{curvatura scalare};
\vspace{0,2cm}
\item $|\sigma_2|^2,\,|\pi_2|^2,$ and a bilinear expression in $\pi_2$, $\sigma_2$ which does not
appear;
\vspace{0,2cm}
\item $|\nu_3|^2$.
\end{itemize}
The 17 copies of $V_{1,1}$ are generated by the projections of
\begin{itemize}
\item $\pi_0\pi_2,\,\pi_0\sigma_2,\,\sigma_0\sigma_2,\sigma_0\pi_2$;
\vspace{0,2cm}
\item 4 bilinear expressions in $\pi_1$  and $\nu_1$ which does not
  appear in formula\\ \eqref{riccitensor};
\vspace{0,2cm}
\item
$*\pi_1\wedge J\nu_3$ and 3 more bilinear expressions in
  $\pi_1$ and $\nu_3$;
\vspace{0,2cm}
\item $*(\pi_2\wedge\pi_2),\,*(\sigma_2\wedge\sigma_2)$ and 2 more bilinear expressions
  in $\pi_2$ and $\sigma_2$;
\vspace{0,2cm}
\item a bilinear form in $\nu_3$.
\end{itemize}
The 12 copies of $V_{2,0}$ are generated by the projections of
\begin{itemize}
\item $\pi_0\nu_3,\,\sigma_0\nu_3$;
\vspace{0,2cm}
\item
  $\nu_1\wedge*(J\pi_1\wedge\Omega),\,J\nu_1\wedge*(\nu_1\wedge\Omega)$
  and other 2 bilinear expressions in
$\pi_1$, $\nu_1$;
\vspace{0,2cm}
\item
  $\sigma_2\wedge\nu_1,\,\pi_2\wedge\nu_1,\,\sigma_2\wedge\pi_1,\,\pi_2\wedge\pi_1$;
\vspace{0,2cm}
\item  two bilinear expressions in $\sigma_2,\,\nu_3$ and
  $\pi_2,\,\nu_3$;
\vspace{0,2cm}
\item $Q(\nu_3,\nu_3)$.
\end{itemize}
An analogous discussion can be done for the second order expressions
after considering the splitting:
$$
\begin{aligned}
V_{2}(\mathfrak{su}(3))=&3\,V_{0,0}\oplus 4\,V_{1,0}\oplus
5\,V_{1,1}\oplus 3\,V_{2,1}\oplus 4\,V_{2,0} \oplus V_{3,0}\oplus V_{2,2ù}\,.
\end{aligned}
$$
}
\end{rem}

\section{The Ricci tensor in the GCY case}
\label{SGCYsection}
Suppose now that the pair $(\kappa,\Omega)$ gives a generalized Calabi-Yau
structure on $M$. In this case all the torsion is encoded by
$\pi_2$ and $\sigma_2$; in fact $d\Omega$ and $dJ\Omega$ reduce to
\begin{equation*}
d\Omega=-\pi_2\wedge\kappa\,,\quad dJ\Omega=-\sigma_2\wedge\kappa\,.
\end{equation*}
Therefore we get
\begin{eqnarray*}
&& 0 = d^2\Omega=-d\pi_2\wedge\kappa\,,\\
&& 0 = d^2J\Omega=-d\sigma_2\wedge\kappa\,,
\end{eqnarray*}
i.e. $d\pi_2$ and $d\sigma_2$ are effective 3-forms.
Since $\pi_2\in \Lambda^2_8M$
$$
\begin{aligned}
0=d(\pi_2\wedge \Omega)&=d\pi_2\wedge \Omega+\pi_2\wedge d\Omega\\
&=d\pi_2\wedge
\Omega-\pi_2\wedge\pi_2\wedge\kappa\\
&=d\pi_2\wedge\Omega+\pi_2\wedge*\pi_2\\
&=d\pi_2\wedge\Omega+|\pi_2|^2*1\,,
\end{aligned}
$$
i.e.
\begin{equation*}
d\pi_2\wedge \Omega=-|\pi_2|^2*1\,.
\end{equation*}
Analogously we get $$d\sigma_2\wedge J\Omega = -|\sigma_2|^2*1\,.$$
Now we can express the Ricci tensor of a generalized Calabi-Yau
manifold in terms of $\pi_2$ and $\sigma_2$. In this case equation
\eqref{riccitensor} reduces to
\begin{equation*}
Ric_0=\frac{1}{4}\,\iota^{-1}(E_1(*(\pi_2 \wedge \pi_2 + \sigma_2\wedge\sigma_2)))
-2\,\gamma^{-1}(E_2(Jd\pi_2+\bigstar d\sigma_2))\,.
\end{equation*}
Since $d\sigma_2$ is effective, $\bigstar d\sigma_2=-d\sigma_2$. Thus
\begin{equation*}
Ric_0=\frac{1}{4}\,\iota^{-1}(E_1(*(\pi_2 \wedge \pi_2 + \sigma_2\wedge\sigma_2)))
-2\,\gamma^{-1}(E_2(Jd\pi_2-d\sigma_2))\,.
\end{equation*}
By the definitions of $E_1$ and $E_2$, using the $J$-invariance of $\pi_2$ and formula
\eqref{ultima}, we have
\begin{equation*}
\begin{aligned}
E_1(*(\pi_2\wedge\pi_2))&=*(\pi_2\wedge\pi_2)
-\frac{1}{9}*((\pi_2\wedge\pi_2
+*(\pi_2\wedge\pi_2)\wedge \kappa)\wedge\kappa)\kappa\\
&=*(\pi_2\wedge\pi_2)+\frac{1}{9}|\pi_2|^2\kappa-
\frac{1}{9}*(*(\pi_2\wedge\pi_2)\wedge\kappa^2)\kappa\\
&=*(\pi_2\wedge\pi_2)+\frac{1}{9}|\pi_2|^2\kappa+\frac{2}{9}|\pi_2|^2\kappa\\
&=*(\pi_2\wedge\pi_2)+\frac{1}{3}|\pi_2|^2\kappa
\end{aligned}
\end{equation*}
and
\begin{equation*}
\begin{aligned}
E_2(d\pi_2)&= d\pi_2-\frac{1}{2}\,*(J d\pi_2\wedge\kappa)\wedge\kappa
                        -\frac{1}{4}*(d\pi_2\wedge J\Omega)\,\Omega
                        +\frac{1}{4}*(d\pi_2 \wedge \Omega)\,J\Omega\\
              &=d\pi_2-\frac{1}{4}\,*(d\pi_2 \wedge J\Omega)\,\Omega
                     -\frac{1}{4}|\pi_2|^2\,J\Omega\\
              &=d\pi_2+\frac{1}{4}\,*(\pi_2 \wedge \sigma_2\wedge
                                      \kappa)\,\Omega
                     -\frac{1}{4}|\pi_2|^2\,J\Omega\,,
\end{aligned}
\end{equation*}
where in the last step we have used 
$$
0=d(\pi_2\wedge J\Omega)=d\pi_2\wedge J\Omega+\pi_2\wedge d J\Omega=
d\pi_2\wedge J\Omega-\pi_2\wedge \sigma_2\wedge\kappa\,.
$$
In the same way we get 
$$
E_1(*(\sigma_2\wedge\sigma_2))=*(\sigma_2\wedge\sigma_2)+\frac{1}{3}|\sigma_2|^2\kappa
$$
and 
$$
E_2(d\sigma_2)=d\sigma_2+\frac{1}{4}*(\pi_2 \wedge \sigma_2\wedge\kappa)\,J\Omega
+\frac{1}{4}|\sigma_2|^2\,\Omega\,.
$$
Therefore, taking into account that $E_2$ commutes with $J$,
the traceless Ricci tensor
of a generalized Calabi-Yau manifold is given
by
\begin{equation}
\label{ricci_0}
\begin{aligned}
Ric_0=&
\frac{1}{4}\,\iota^{-1}(*(\sigma_2\wedge\sigma_2+\pi_2\wedge\pi_2)+
\frac{1}{3}(|\sigma_2|^2+|\pi_2|^2)\,\kappa)\\
&-2\gamma^{-1}(Jd\pi_2-d\sigma_2+\frac{1}{4}(|\pi_2|^2-|\sigma_2|^2)\,\Omega)\,.
\end{aligned}
\end{equation}
Formula \eqref{ricci_0} implies that the metric induced by a GCY structure $(\kappa,\Omega)$
is Einstein ({\it i.e.} $Ric_0=0$) if and only if the torsion forms $\pi_2,\sigma_2$ satisfies
\begin{equation}
\label{drag}
\begin{cases}
\sigma_2\wedge\sigma_2+\pi_2\wedge\pi_2+
\frac{1}{6}(|\pi_2|^2+|\sigma_2|^2)\,\kappa\wedge\kappa=0\\[3pt]
Jd\pi_2-d\sigma_2+\frac{1}{4}(|\pi_2|^2-|\sigma_2|^2)\,\Omega=0\,.
\end{cases}
\end{equation}

In the special case of SGCY manifolds we can prove
\begin{cor}
\label{Cleyton}
A $6$-dimensional {\em SGCY} manifold is Einstein if and only if it is a genuine Calabi-Yau 
manifold.
\end{cor}

The proof of Corollary \ref{Cleyton} relies on the following lemma which is 
interesting in its own.


\begin{lemma}\label{Richard}
Let $(V,\kappa,\Omega)$ be a $6$-dimensional symplectic vector space endowed with a normalized 
$\kappa$-positive $3$-form. If $\alpha\neq 0$ belongs to $\Lambda^2_8V^*$, then 
$\alpha \wedge \alpha$ does not belong to the $1$-dimensional {\em SU(3)}-module generated by 
$\kappa \wedge \kappa$.
\end{lemma} 
\begin{proof}
The key observation here is that $\Lambda^2_8 V^*$ is isomorphic as a SU$(3)$-representation 
to the adjoint representation $V_{1,1}$. Since every element in $\mathfrak{su}(3)$ is 
$\mbox{Ad}(\mbox{SU}(3))$-conjugated to an element of a fixed Cartan subalgebra of 
$\mathfrak{su}(3)$, there exists a SU$(3)$-basis $\{e^1,\ldots,e^6\}$ of $V^*$ such that
\[
\begin{aligned}
&\alpha=\lambda_1\,e^{12}+\lambda_2\,e^{34}-(\lambda_1+\lambda_2)\,e^{56}\,,
\end{aligned}
\]   
for some $\lambda_1,\lambda_2 \in \R$. Now suppose that 
$\alpha\wedge\alpha = q\,\kappa\wedge\kappa$ for some $q \in \R$.
Setting to zero the three components of $\alpha \wedge \alpha -q\,\kappa\wedge\kappa$ gives 
the equations
\begin{eqnarray*}
&&\lambda_1^2+\lambda_1\lambda_2+q=0\,,\\
&&\lambda_2^2+\lambda_1\lambda_2+q=0\,,\\
&&\lambda_1\lambda_2-q=0\,,
\end{eqnarray*}
which readily imply $q=0$.
\end{proof}
\begin{proof}[Proof of corollary \ref{Cleyton}]
Since in the GCY case $\pi_2=0$, taking into account lemma \ref{Richard}, 
the first equation of \eqref{drag} can be satisfied if and only if 
$|\sigma_2|^2=0$. Therefore the Einstein condition forces $(\kappa,\Omega)$ to be a 
Calabi-Yau structure on $M$.  
\end{proof}
\begin{rem}
\emph{In \cite{Drag} it has been proven (see theorem 1) that a} compact 
\emph{Einstein almost
K\"ahler manifold with vanishing first Chern class
is actually a K\"ahler-Einstein manifold. Note that our result holds with no the compactness 
assumption.}
\end{rem}

\section{An explicit example}
\label{examplesection}
In this last section we carry out the computation of the Ricci tensor
and the intrinsic torsion
of a left-invariant SU$(3)$-structure on a particular 6-dimensional
nilmanifold.
\newline

Let G be the nilpotent Lie group of the matrices of the form
$$
A=\left(
\begin{array}{cccccc}
1&0&x_1&x_3&0&0\\
0&1&x_2&x_4&0&0\\
0&0&1&x_5&0&0\\
0&0&0&1&0&0\\
0&0&0&0&1&x_6\\
0&0&0&0&0&1
\end{array}
\right)\,
$$
where $x_1,x_2,x_3,x_4,x_5,x_6$ are real numbers. Let $\Gamma$ be the set
of matrices in G having integral entries, then $M:=$G$/\Gamma$ is a
compact parallelizable  smooth manifold. Let
$\{X_1,\dots,X_n\}$ be the global  frame on $M$ given
 by
\begin{equation*}
\begin{aligned}
&X_1=\de{}{x_5}+x_1\de{}{x_3}+x_2\de{}{x_4}\,,\;\;X_2=
\de{}{x_6}\,,\\
&X_3=\de{}{x_2}\,,\;\;
X_4=\de{}{x_3}\,,\;\;X_5=\de{}{x_1}\,,\;\;X_6=\de{}{x_4}\,.
\end{aligned}
\end{equation*}
We have that
$$
[X_1,X_3]=-X_6\,,\quad [X_1,X_5]=-X_4
$$
and the other brackets are zero. Let $\{\alpha_1,\dots,\alpha_6\}$ be
the dual frame of $\{X_1,\dots,X_n\}$, then
$$
\begin{cases}
d\alpha_1=d\alpha_2=d\alpha_3=d\alpha_5=0\\
d\alpha_4=\alpha_{15}\\
d\alpha_6=\alpha_{13}\,.
\end{cases}
$$
Therefore the \emph{closed}  global forms
$$
\begin{aligned}
&\kappa=\alpha_{12}+\alpha_{34}+\alpha_{56}\,,\\
&\Omega=\alpha_{135}-\alpha_{146}-\alpha_{245}-\alpha_{236}\,.
\end{aligned}
$$
defines a SGCY structure on $M$. Let $J$ be the almost complex
structure on $M$ induced by the SU(3)-structure, then on the frame $\{X_1,\dots,X_6\}$ one
has
$$
J(X_1)=X_2\,,\quad J(X_{3})=X_4\,,\quad J(X_5)=X_6\,.
$$
We have
$$
dJ\Omega=d(-\alpha_{246}+\alpha_{235}+\alpha_{145}+\alpha_{136})=\alpha_{1234}-\alpha_{1256}=(\alpha_{34}-\alpha_{56}
)\wedge\kappa\,,
$$
i.e., with the notations of \eqref{scomposition},
$$
\sigma_2=\alpha_{56}-\alpha_{34}\,.
$$
Since $(M,\kappa,\Omega)$ is a SGCY manifold, $\sigma_2$ is the only
non-zero torsion form.\\
Note that  the metric associated to
$(\kappa,\Omega)$ is
$$
g=\sum_{i=1}^n \alpha_i\otimes\alpha_i\,.
$$
Consequently we have $|\sigma_2|^2=2$, hence formula \eqref{curvatura scalare SGCY}
 implies $s=-1$.

Using \eqref{ricci_0} we can compute the Ricci tensor of $g$: we have
$$
\begin{aligned}
Ric_0=&\iota^{-1}(-\frac{1}{2}\,\alpha_{12}+\frac{1}{6}\,\kappa)+
\gamma^{-1}(-4\,\alpha_{135}+\Omega)\\
=&\iota^{-1}(-\frac{1}{3}\,\alpha_{12}+\frac{1}{6}\,\alpha_{34}+\frac{1}{6}\,\alpha_{56})+\\
&+\gamma^{-1}(-3\,\alpha_{135}-\,\alpha_{146}-\,\alpha_{245}-\,\alpha_{236})\,.\\
\end{aligned}
$$
Let $\nabla$ be the Levi-Civita connection of $g$, then
\begin{equation*}
\begin{aligned}
&\nabla_{1}X_3=-\frac{1}{2}X_6\,,&&\nabla_{1}X_6=\frac{1}{2}X_3\,,&&&\nabla_{3}X_6=-\frac{1}{2}X_1\,,\\
&\nabla_{3}X_1=\frac{1}{2}X_6\,,&&\nabla_{6}X_1=\frac{1}{2}X_3\,,&&&\nabla_{6}X_3=-\frac{1}{2}X_1\,,\\
&\nabla_{1}X_5=-\frac{1}{2}X_4\,,&&\nabla_{1}X_4=\frac{1}{2}X_5\,,&&&\nabla_{5}X_4=-\frac{1}{2}X_1\,,\\
&\nabla_{5}X_1=\frac{1}{2}X_4\,,&&\nabla_{4}X_1=\frac{1}{2}X_5\,,&&&\nabla_{4}X_5=-\frac{1}{2}X_1\,,
\end{aligned}
\end{equation*}
where $\nabla_{i}X_j$ stands for $\nabla_{X_i}X_j$. Now are ready to
compute the torsion of this SU(3)-manifold. We immediately have
$$
\psi=\frac{1}{2}\,\left(
\begin{array}{cccccc}
0&0&-\alpha_6&-\alpha_5&-\alpha_4&-\alpha_3\\
0&0&0&0&0&0\\
\alpha_6&0&0&0&0&\alpha_1\\
\alpha_5&0&0&0&-\alpha_1&0\\
\alpha_4&0&0&\alpha_1&0&0\\
\alpha_3&0&-\alpha_1&0&0&0
\end{array}
\right)
$$
and a computation gives
$$
\theta=\frac{1}{4}\left(
\begin{array}{cccccc}
0&0&-\alpha_6&-\alpha_5&-\alpha_4&-\alpha_3\\
0&0&\alpha_5&-\alpha_6&\alpha_3&-\alpha_4\\
\alpha_6&-\alpha_5&0&0&0&2\alpha_1\\
\alpha_5&\alpha_6&0&0&-2\alpha_1&0\\
\alpha_4&-\alpha_3&0&2\alpha_1&0&0\\
\alpha_3&\alpha_4&-2\alpha_1&0&0&0
\end{array}
\right)
$$
and
$$
\tau=\frac{1}{4}\left(
\begin{array}{c}
0\\
0\\
\alpha_5\\
-\alpha_3\\
-\alpha_6\\
\alpha_5
\end{array}\right)
\,,\quad \mu=0\,.
$$
\section{Appendix}
\newcommand\si{*_{\sigma}}
In this appendix we give a proof of lemma \ref{Bryantlemma} and theorem \ref{t s}.
\begin{proof}[Proof of lemma \ref{Bryantlemma}]
Let $N$ be the Riemannian product $N=M\times\R$. Denote by
$$
\begin{aligned}
&p_1\colon N\to M\,,\\
&p_2\colon N\to \R\\
\end{aligned}
$$
the projections. The 3-form
$$
\sigma=p_1^*(\Omega)+p_1^*(\kappa)\wedge p_2^*(dt)\,,
$$
defines a G$_2$-structure on $N$. From now on we identify the forms
$\kappa$, $\Omega$, $dt$ with their respective pull-backs to $N$. Let us
denote by $\si$ and $*$ the Hodge operator associated to the metric induced
by $\sigma$ and by the SU$(3)$-structure on $M$ respectively.
Thus
\begin{eqnarray*}
&&d\sigma=d\Omega+d\kappa\wedge dt\,,\\
&&\si\s=(*\Omega)\wedge dt+\,* \kappa=J\Omega\wedge dt+\frac{1}{2}\,\kappa^2\,,\\
&&d\si\s=dJ\Omega\wedge dt+\,d\kappa\wedge\kappa\,,\\
&&\si d\s=(*d\Omega)\wedge dt-\,* d\kappa\,,\\
&&\si d\si\sigma=*dJ\Omega+\,*(d\kappa\wedge\kappa)\wedge dt\,.
\end{eqnarray*}
Now we use the formula
\begin{equation}
\label{bryant}
\si\s\wedge\si(d\si\s)+(\si d\s)\wedge\s=0\,,
\end{equation}
proved by Bryant in \cite{B1}. Now we have
$$
\begin{aligned}
&\si\s\wedge\si(d\si\s)+(\si d\s)\wedge\s=J\Omega\wedge(*dJ\Omega)\wedge
dt+\frac{1}{2}\,\kappa^2\wedge(*(d\kappa\wedge\kappa))\wedge dt+\\
&+\frac{1}{2}\,\kappa^2\wedge*dJ\Omega-\,(*d\Omega)\wedge\Omega\wedge dt
-\,(*d\kappa)\wedge\Omega-\,(*d\kappa)\wedge\kappa\wedge dt \,.
\end{aligned}
$$
Therefore equation \eqref{bryant} implies
\begin{itemize}
\item $(*d\kappa)\wedge\Omega=\frac{1}{2}\,\kappa^2\wedge*dJ\Omega$, which is indeed
an easy consequence of $\Omega \wedge \kappa=0$:
\vspace{0.1 cm}
\item $J\Omega\wedge(*dJ\Omega)+\frac{1}{2}\,\kappa^2\wedge*(d\kappa\wedge\kappa)
-\,(*d\Omega)\wedge\Omega
-\,(*d\kappa)\wedge\kappa =0\,.
$
\end{itemize}
In order to show that equation
\eqref{Bryant 2}
holds, we need to prove the following identity
\begin{equation}
\label{identity}
\frac{1}{2}\,\kappa^2\wedge*(d\kappa\wedge\kappa)=(*d\kappa)\wedge\kappa \,.
\end{equation}
The decomposition of 3-forms on $M$ implies
\begin{equation*}
\frac{1}{2}\,\kappa^2\wedge*(d\kappa\wedge\kappa)=\frac{1}{2}\,\kappa^2\wedge*(\nu_1\wedge\kappa^2)=(\bigstar\kappa)
\wedge*(\nu_1\wedge\kappa^2)
\end{equation*}
and
\begin{equation*}
(*d\kappa)\wedge\kappa=(*(\nu_1\wedge\kappa))\wedge\kappa\,,
\end{equation*}
where
$\nu_1\wedge\kappa\in\Lambda_{6}^3M=\{\gamma\in \Lambda^3M\,|\,\bigstar\gamma=\gamma\}$.
Now we need to recall the following lemma proved in \cite{dBT1};\\[7pt]
\textbf{Lemma A.1.}\emph{
Let $\zeta\in \Lambda^1V^*$ and $\gamma\in \Lambda^rV^*$; we have
\begin{equation}
\label{Db-T}
\bigstar(\zeta\wedge\gamma)=(-1)^r\zeta\wedge\bigstar(\kappa\wedge\gamma)-(-1)^r\bigstar(\kappa\wedge\bigstar
(\zeta\wedge\bigstar\gamma))\,.
\end{equation}}
\vspace{0.3 cm}
Applying equation \eqref{Db-T} with $\zeta=*(\nu_1\wedge\kappa^2)$ and
$\gamma=1\in\Lambda^0M$ we have
\begin{equation}
\label{primai}
(\bigstar\kappa)\wedge*(\nu_1\wedge\kappa^2)=\bigstar(*(\nu_1\wedge\kappa^2))=*J(*(\nu_1\wedge\kappa^2))=
-\,J\nu_1\wedge\kappa^2\,.
\end{equation}
Moreover, since $\nu_1\in\Lambda^3_6M$, it follows
\begin{equation}
\label{secondai}
*(\nu_1\wedge\kappa)\wedge\kappa=-J\nu_1\wedge\kappa^2\,.
\end{equation}
Equation \eqref{primai} together with equation \eqref{secondai}  imply
\eqref{identity}, so that equation \eqref{Bryant 2} is proved.
\end{proof}
\begin{proof}[Proof of theorem \ref{t s}.]
In order to prove formula \eqref{curvatura scalare} it is useful to
introduce the 1-forms $S_{ijk}\,\omega_k$, $V_{ik}\,\omega_k$,
defined by the relations
\begin{eqnarray*}
&&dT_{ij}=T_{ik}\,\theta_{kj}+T_{kj}\,\theta_{ki}+S_{ijk}\,\omega_k\,,\\
&&dM_i=M_k\,\theta_{ki}+V_{ik}\,\omega_{k}\,.
\end{eqnarray*}
Using equations \eqref{Dtau} and \eqref{Dmu} and the definition of
$T_{ij}$, $M_{i}$ given in \eqref{deftaumu}
$$
\begin{aligned}
D\tau_i=&dT_{ij}\wedge\omega_{j}+T_{ij}\,d\omega_j-2\,\kappa_{ij}\,\mu\wedge\tau_j\\
=&(S_{iba}-T_{ij}T_{qa}\e_{jbq}-T_{ij}\kappa_{jb}M_a-2\kappa_{ij}M_aT_{jb})\,\omega_a\wedge\omega_b\,,
\end{aligned}
$$
and
$$
\begin{aligned}
D\mu=&dM_r\wedge\omega_r+M_rd\omega_r+\frac{2}{3}\kappa_{ij}\,\tau_{i}\wedge\tau_j\\
=&(V_{ba}-M_r\e_{rbq}T_{qa}-M_r\kappa_{rb}M_a+\frac{2}{3}\kappa_{ij}T_{ia}T_{jb})\,\omega_a\wedge\omega_b\,.
\end{aligned}
$$
Therefore, taking into account \eqref{full}, \eqref{coppia}, we obtain
$$
\begin{aligned}
T_{iab}&=2(S_{iba}-T_{ij}T_{qa}\e_{jbq}-T_{ij}\kappa_{jb}M_a-2\kappa_{ij}M_aT_{jb})\,,\\
N_{ab}&=2(V_{ba}-M_r\e_{rbq}T_{qa}-M_r\kappa_{rb}M_a+\frac{2}{3}\kappa_{ij}T_{ia}T_{jb})\,.
\end{aligned}
$$
It follows that
\begin{eqnarray*}
&&\e_{ipq}T_{pqj}=2(\e_{ipq}S_{pjq}-
\e_{ipq}\e_{rjs}T_{pr}T_{sq}-\e_{ipq}T_{pr}\kappa_{rj}M_{q}+2\ov{\e}_{iqr}T_{rj}M_q)\,,\\
&&\kappa_{ip}N_{pj}=2(\kappa_{ip}V_{jp}-\kappa_{ip}\e_{rjq}T_{qp}M_r-\kappa_{ip}\kappa_{rj}M_{r}M_p
+\frac{2}{3}\kappa_{ip}\kappa_{qr}T_{qp}T_{rj})\,
\end{eqnarray*}
and using the $\epsilon$-identities \eqref{ep}
\begin{equation*}
\begin{aligned}
\e_{ipq}T_{pqi}=&2(-\e_{ipq}S_{ipq}-
\e_{ipq}\e_{ris}T_{pr}T_{sq}-\ov{\e}_{prq}T_{pr}M_{q}+2\ov{\e}_{qri}T_{ri}M_q)\\[3pt]
=&2(-\e_{ipq}S_{ipq}-
\e_{ipq}\e_{ris}T_{pr}T_{sq}+\ov{\e}_{prq}T_{pr}M_{q})\,,\\[3pt]
\kappa_{ip}N_{pi}=&2(\kappa_{ip}V_{ip}-\kappa_{ip}\e_{riq}T_{qp}M_r-\kappa_{ip}\kappa_{ri}M_{r}M_p
+\frac{2}{3}\kappa_{ip}\kappa_{qr}T_{qp}T_{ri})\\[3pt]
=&2(\kappa_{ip}V_{ip}+\ov{\e}_{rqp}T_{qp}M_r
+\frac{2}{3}\kappa_{ip}\kappa_{qr}T_{qp}T_{ri}+\Sigma_i\, M_i^2) \,.
\end{aligned}
\end{equation*}
Then by theorem \ref{controbuio} we get
\begin{equation*}
\begin{aligned}
s=&4(-\e_{ipq}S_{ipq}-
\e_{ipq}\e_{ris}T_{pr}T_{sq}+\ov{\e}_{prq}T_{pr}M_{q})\\[3pt]
&-6(\kappa_{ip}V_{ip}+\ov{\e}_{rqp}T_{qp}M_r
+\frac{2}{3}\kappa_{ip}\kappa_{qr}T_{qp}T_{ri}+\Sigma_i\, M_i^2)\\[3pt]
=&-4\e_{ipq}S_{ipq}-
4\e_{ipq}\e_{ris}T_{pr}T_{sq}-2\ov{\e}_{prq}T_{pr}M_{q}\\[3pt]
&-6\kappa_{ip}V_{ip}
-4\kappa_{ip}\kappa_{qr}T_{qp}T_{ri}-6\Sigma_i\, M_i^2\,.
\end{aligned}
\end{equation*}
Furthermore a straightforward computation  gives the following
formulae
\begin{equation*}
\begin{aligned}
&\pi_0^2=\frac{4}{9}T_{ii}T_{jj}\,,\\[3pt]
&\sigma_0^2=\frac{4}{9}\kappa_{ij}\kappa_{sr}T_{ij}T_{sr}\,,\\[3pt]
&|\pi_2|^2=-\frac{4}{3}T_{ii}T_{jj}+4T_{ij}^2-2\e_{sra}\e_{aij}T_{sr}T_{ij}+4\kappa_{ir}\kappa_{js}
T_{ij}T_{sr}\,,\\[3pt]
&|\sigma_2|^2=-2\e_{sra}\e_{aij}T_{sr}T_{ij}
-\frac{4}{3}\kappa_{ij}\kappa_{ab}T_{ij}T_{ab}-4T_{ij}T_{ji}+4\Sigma_{ij}T_{ij}^2\,,\\[3pt]
&|\nu_1|^2=\e_{ijk}\e_{kab}T_{ij}T_{ab}\,,\\[3pt]
&|\nu_3|^2=2T_{ij}^2+2T_{ij}T_{ji}-2\kappa_{jr}\kappa_{is}T_{ij}T_{rs}-2\kappa_{ir}\kappa_{js}T_{ij}T_{rs}\,,\\[3pt]
&d^*\pi_1=-\e_{sra}\e_{aij}T_{sr}T_{ij}+4\ov{\e}_{ijk}T_{ij}M_k-\e_{sra}S_{sra}-3\kappa_{ij}V_{ij}-3\Sigma_i\,M_i^2\,,\\[3pt]
&d^*\nu_1=-\e_{sra}\e_{aij}T_{sr}T_{ij}+\ov{\e}_{ijk}T_{ij}M_k-\e_{sra}S_{sra}\,,\\[3pt]
&\langle\pi_1,\nu_1\rangle=\e_{abk}\e_{kij}T_{ab}T_{ij}-3\ov{\e}_{ijk}T_{ij}M_k\,.
\end{aligned}
\end{equation*}
Therefore we get
$$
\begin{aligned}
&\frac{15}{2}\pi_0^2+\frac{15}{2}\sigma_0^2+2d^*\pi_1+2d^*\nu_1
-|\nu_1|^2-\frac{1}{2}|\sigma_2|^2
-\frac{1}{2}|\pi_2|^2
-\frac{1}{2}|\nu_3|^2+4\langle\pi_1,\nu_1\rangle=\\[3pt]
=&4T_{ii}T_{jj}+4\kappa_{ij}\kappa_{sr}T_{ij}T_{sr}-5\Sigma_{ij}T_{ij}+\e_{sra}\e_{aij}T_{sr}T_{ij}
+T_{ij}T_{ji}-2\ov{\e}_{ijk}T_{ij}M_k\\[3pt]
&-6\kappa_{ij}V_{ij}-6\Sigma_i\,M_i^2+(-\kappa_{ia}\kappa_{jb}+\kappa_{ib}\kappa_{ja})T_{ij}T_{ba}-4\e_{ijk}S_{ijk}=\\
=&4\e_{ipq}S_{ipq}-
4\e_{ipq}\e_{ris}T_{pr}T_{sq}-2\ov{\e}_{prq}T_{pr}M_{q}
-6\kappa_{ip}V_{ip}
-4\kappa_{ip}\kappa_{qr}T_{qp}T_{ri}-6\Sigma_i\, M_i^2\,,
\end{aligned}
$$
i.e.
$$
s=\frac{15}{2}\pi_0^2+\frac{15}{2}\sigma_0^2+2d^*\pi_1+2d^*\nu_1
-|\nu_1|^2-\frac{1}{2}|\sigma_2|^2
-\frac{1}{2}|\pi_2|^2
-\frac{1}{2}|\nu_3|^2+4\langle\pi_1,\nu_1\rangle
\,,
$$
and the  theorem is proved.
\end{proof}

\end{document}